\documentclass[11pt]{article}

\usepackage{algpseudocode,lineno}
\usepackage{enumitem}
\usepackage{amssymb,multirow,color,comment,amsmath}
\usepackage{amsfonts,amsthm, bm}
\usepackage{xcolor}
\usepackage[font=itshape]{quoting}
\usepackage{algorithm,algpseudocode}
\usepackage{verbatim}
\usepackage{amssymb,multirow,color,comment,amsmath}
\usepackage{amsthm,tikz}
\usetikzlibrary{matrix}
\usepackage{pgfplots,cite}
\usepackage{subfigure,bm}
\usepackage{algorithm,algpseudocode}
\usepackage{breqn}
\usepackage[font=itshape]{quoting}

\usepackage[margin=1.4in]{geometry}

\usetikzlibrary{arrows,calc}
\usepackage{multirow,comment}
\usepackage{pgf,tikz,tikz-3dplot}
\usepackage[percent]{overpic}
\usepackage{pdflscape}
\usepackage{pgfplots,xcolor}
\usepackage[makeroom]{cancel}
\numberwithin{equation}{section}

\newcommand{\br}{\mathbf{r}}
\newcommand{\sd}{\mbox{d}}

\newcommand{\bx}{\mathbf{x}}

\newcommand{\by}{\mathbf{y}}
\newcommand{\red}[1]{{\color{black} #1}}
\DeclareMathAlphabet\mathbfcal{OMS}{cmsy}{b}{n}

\begin{document}

\title{\red{Accelerating frequency-domain numerical methods for
weakly nonlinear focused ultrasound using nested meshes}}
\author{Samuel P.~Groth$^1$, Pierre G\'{e}lat$^2$, Seyyed R.~Haqshenas$^{2}$, \\
Nader Saffari$^2$, Elwin van 't Wout$^3$, Timo Betcke$^4$, Garth N.~Wells$^1$\\[6pt]
$^1$Department of Engineering, University of Cambridge \\[4pt]
$^2$Department of Mechanical Engineering, University College London, \\[4pt]
$^3$Institute for Mathematical and Computational Engineering, \\
School of Engineering and Faculty of Mathematics, \\
            Pontificia Universidad Cat\'olica de Chile \\ [4pt]
$^4$Department of Mathematics, University College London}

\date{}
\maketitle
\begin{abstract}
    The numerical simulation of \red{weakly} nonlinear ultrasound is
    important in treatment planning for \red{focused ultrasound (FUS)}
    therapies. However, the large domain sizes and generation of higher
    harmonics at the focus make these problems extremely computationally
    demanding. Numerical methods typically employ a uniform mesh fine
    enough to resolve the highest harmonic present in the problem,
    leading to a very large number of degrees of freedom. This paper
    proposes a more efficient strategy in which each harmonic is
    approximated on a separate mesh, the size of which is proportional
    to \red{the} wavelength of the harmonic. The increase in resolution
    required to resolve a smaller wavelength is balanced by a reduction
    in the domain size. This nested meshing is feasible owing to the
    increasingly localised nature of higher harmonics near the focus.

    Numerical experiments are performed for \red{FUS} transducers in
    homogeneous media in order to determine the size of the meshes
    required to accurately represent the harmonics. In particular, a
    fast \textit{volume potential} approach is proposed and employed to
    perform convergence experiments as the computation domain size is
    modified. This approach allows each harmonic to be computed via the
    evaluation of an integral over the domain. Discretising this
    integral using the midpoint rule allows the computations to be
    performed rapidly with the FFT. It is shown that at least an order
    of magnitude reduction in memory consumption and computation time
    can be achieved with nested meshing. \red{Finally, it is
    demonstrated how to generalise this approach to inhomogeneous
    propagation domains.}
\end{abstract}

\clearpage
\section{\label{sec:1} Introduction}

\red{Focused ultrasound (FUS)} is a non-invasive tumor ablation therapy
in which acoustic waves are focused at a target location, thereby
elevating the temperature sufficiently to destroy the tumor tissue.
Tissue death can be caused directly via thermal ablation, or via other
mechanisms such as
cavitation~\cite{vlaisavljevich2013image,izadifar2017mechanical},
\red{shock-scattering histotripsy~\cite{maxwell2011cavitation}, and
boiling histotripsy~\cite{khokhlova2011controlled}}. In the thermal
ablation setting, which is the focus of this article, the peak acoustic
pressure is often between 1 and 10 MPa, at which linear acoustic theory
is no longer accurate. Indeed it has been observed that the
contributions from nonlinear effects can increase the temperature at the
focal point by an additional 20\%, thus substantially reducing the time
required for ablation when compared to simulations using linear
theory~\cite{solovchuk2014temperature}.
To simulate the nonlinear acoustic field, a popular model is the
Westervelt equation, which incoporates a quadratic pressure term.
Numerically solving the Westervelt equation in a \red{FUS} setting is
computationally challenging since the domains are large compared to the
smallest wavelength present.

To put into context the scale of the computational problems presented by
\red{FUS}, consider a bowl-shaped \red{FUS} transducer with radius 7.5cm
and geometric focal length 15cm operating at 1.1MHz. Transducers of this
size are suitable for deep-seated tumors, located in the liver, for
example. A simulation domain spanning the face of the bowl and extending
to the focal region is approximately 100 wavelengths across in each
dimension. If there are five harmonics present in the nonlinear field,
then this is five hundred wavelengths at the fifth harmonic. In order to
obtain a reasonable accuracy level, we can assume that at least six
degrees of freedom (DOF) per wavelength are
required~\cite{marburg2008discretization}. Then to simulate this problem
in three dimensions, we would require 27 billion degrees of freedom.
Handling such a large system presents an enormous computational load in
terms of memory and time.

There has been a great deal of research effort aimed at reducing the
computational cost of nonlinear \red{FUS} simulations. Much of this work
employs at least one simplifying assumption, such as axisymmetry,
one-way wave propagation, or the parabolic approximation (see
\cite{gu2015modeling} for a review). Popular methods for solving the
original Westervelt equation (or equivalent) include the
finite-difference time-domain method~\cite{solovchuk2013simulation}
(FDTD) and the $k$-space pseudospectral method~\cite{treeby2012modeling}
(KSPS). Although powerful techniques, these full-wave solvers can still
yield simulation times in excess of a day for realistic problems when
run on a large cluster~\cite{jaros2016full}. \red{For efficiently
simulating highly-focused problems with strong nonlinearities,
shock-capturing numerical schemes have been developed and applied for
the parabolic approximation of the Westervelt equation, known as the KZK
equation~\cite{bessonova2009focusing}, wide-angle parabolic
approximation~\cite{yuldashev2018wide}, as well as for the one-way
propagation version of the Westervelt
equation~\cite{yuldashev2011simulation}.}

In this paper, we set out to improve upon the meshing strategies
employed in popular full-wave solvers for the Westervelt equation.
\red{These solvers typically use} uniform meshes resolved according to
the wavelength of the highest harmonic. Such a meshing approach is
inefficient since the higher harmonics are localised around the focus,
therefore the extremely fine mesh far away from the focus is likely
overkill. It seems intuitive that a mesh that is gradually refined as
the focus is approached would be sensible, but at what rate? And what
saving can one expect to achieve? Here we address these questions via
numerical experimentation on realistic \red{FUS} transducers taken from
\cite{sonic}.

\red{We note similar strategies have previously been employed to reduce
the computational load of nonlinear acoustics
simulations~\cite{yuldashev2011simulation,karzova2017shock}. In these
prior works, it was also observed that far from the transducer focus,
higher harmonics provide negligible contributions to the acoustic field
and therefore, to save on storage requirements, higher harmonics are
stored only on smaller domains, centred at the focus. However, the grid
step used in these works is not altered for the different harmonic
computations. By increasing the grid step for lower harmonics, further
computational savings can be achieved, as we demonstrate in this
article. Furthermore, we detail an in-depth study of how to choose
appropriate computational domains for each harmonic and we provide a
useful rule of thumb relating the domain dimensions to each harmonic's
wavelength.}

In order to perform these experiments efficiently, we consider a
simplified setting in which the propagation medium is homogeneous and
the signal is assumed to have harmonic time dependence.
\red{Furthermore, we consider settings in which the peak amplitude is no
more than 15MPa and thus the field is only weakly non-linear}. Under the
time-harmonic \red{and weakly non-linear assumptions}, the full-wave
Westervelt equation reduces to a series of inhomogeneous Helmholtz
equations, one for each harmonic present in the field (as in, e.g.,
\cite{du2013fast}). \red{We note that the weak non-linearity permits us
to reasonably neglect the transfer of energy from higher to lower
harmonics. This is, however, not valid for extremely high amplitude
fields, as encountered in histotripsy applications.}

Since we consider a homogeneous medium, each Helmholtz equation is
exactly solved by the evaluation of the appropriate \textit{volume
potential} integral~\cite{costabel2015spectrum}. Therefore, for each
harmonic in the field, we look at the convergence of the volume
potential as the domain of integration is increased. In order to
efficiently evaluate the volume potential, we discretise the domain over
a uniform voxel mesh, which allows the summation (i.e., the discrete
form of the integration) to be performed rapidly using the fast-Fourier
transform (FFT). This is an extremely efficient technique for computing
high-order harmonics in a homogeneous medium and its application in this
area is, to the best of the authors' knowledge, novel.

For numerical investigations with different transducer configurations
and within two different media (water and liver), we deduce that (for
the configurations considered) in order to perform accurate computations
of the second harmonic, the computation domain must extend all the way
back to the transducer, however may be contracted slightly in the
transverse direction. For the third and higher harmonics, the
computation domain may be contracted in both the axial and transverse
directions, thus leading to substantial computational savings. In fact,
we demonstrate that an accurate approximation (less than 1\% relative
error) may be obtained while contracting the width and length of the
domain in proportion to the wavelength of the harmonic. Thus the number
of cells in every mesh is roughly equal. This amounts to a reduction in
the number of DOF by a factor of approximately $(n/2)^3$, where $n$ is
the number of harmonics being computed.

In practice, this leads to a series of nested meshes, each at the
resolution required for the appropriate harmonic and all with the same
number of voxels. To perform computations for higher harmonics,
solutions from lower harmonics are interpolated onto the finer meshes.
We outline an algorithm for this more efficient computation of all the
harmonics in Section~\ref{subsec:interpolation} and examine its
performance.

The layout of the paper is as follows.
Section~\ref{sec:nonlinear_acoustics} outlines the mathematical model we
employ for our \red{FUS} setup. In particular we consider the Westervelt
equation and review how it reduces to a series of Helmholtz equations
under the time-harmonic assumption. We further simplify the equations
via an assumption of weak nonlinearity. In Section~\ref{sec:volume} we
describe how, in the homogeneous domain case, the Helmholtz equations
are solved exactly via volume potentials. It is then described how these
volume potentials are efficiently approximated using a voxelised
discretisation approach. This leads to the discrete versions of the
volume potentials having block-Toeplitz form; thus the potentials may
each be evaluated using a single fast matrix-vector product with the
FFT. Section~\ref{subsec:incident} discusses our model for the
time-harmonic bowl-shaped transducer. Section~\ref{sec:validation}
presents a validation of our approach by comparing to simulations
performed with the HITU Simulator \textsc{Matlab}
toolbox~\cite{HITUwebpage,soneson2017extending}. In
Section~\ref{sec:results} we perform convergence tests for each of the
harmonics for a range of problem setups. We present our findings on the
rates of convergence of the approximation as the domain is increased and
present a rule of thumb for designing the meshes to ensure each harmonic
is accurately approximated. In Section~\ref{subsec:interpolation} we
present how the hierarchy of meshes is used in practice with our volume
potential approach, including the interpolation of approximations
between the meshes. Here we present some performance details for this
three-dimension full-wave approach on a single workstation. Finally, in
Section~\ref{sec:conclusion} we present our conclusions and discuss the
relevance of this work to other numerical techniques for \red{FUS}.

\section{\label{sec:2} Nonlinear acoustics in the frequency domain}
\label{sec:nonlinear_acoustics}

Acoustic fields produced by \red{FUS} transducers are commonly modelled
by the Westervelt equation~\cite{hamilton1998nonlinear}:
    \begin{equation}
        \nabla^2 p
        - \frac{1}{c_0^2}\frac{\partial^2 p}{\partial t^2}
        - \frac{2\alpha_0}{c_0^{1-\eta}}\frac{\partial }{\partial t}
        \left( -\nabla^2\right)^{\eta/2}p =
        -\frac{\beta}{\rho_0 c_0^4}\frac{\partial^2 p^2}{\partial t^2},
        \label{eqn:west_fractional}
    \end{equation}
    where $p$ is the acoustic pressure, $c_0$ is the speed of sound,
    $\rho_0$ is the medium density, $\beta$ the non-linearity parameter, and
    $\alpha_0$ and $\eta$ are medium specific attenuation parameters.
    The fractional Laplacian appearing in (\ref{eqn:west_fractional}) was first
    introduced in \cite{chen2004fractional} in order to incorporate frequency
    dependent power law attenuation. More specifically, in the frequency domain, it
    generates a complex wavenumber of the form
    \begin{equation}
        k = \frac{\omega}{c_0} + \text{i}\alpha;\quad  \alpha=\alpha_0|\omega|^{\eta},
        \label{eqn:power_law}
    \end{equation}
    for $\alpha_0|\omega|^{\eta-1}c_0<0.1$ ~\cite{szabo1994time}, where
    $\omega$ is the angular frequency of the transducer. The power law
    exponent is typically in the range $1\leq \eta\leq 2$. We note that
    this power law attenuation can also be incorporated via a temporal
    convolution, as was originally proposed in \cite{szabo1994time}; we
    refer the reader to \cite{treeby2010modeling} for a review of power
    law attenuation techniques. \red{We note that the relation
    (\ref{eqn:power_law}) does not incorporate the effect of dispersion
    on the real part of the wavenumber. However, for the frequencies and
    materials considered in this article, this effect is negligibly
    small. For highly nonlinear and higher frequency problems, the
    relation may be modified to more accurately incorporate dispersion,
    as discussed in \cite{treeby2010modeling,waters2005causality}.}

In this article, we assume that the operation time of the transducer is
long when compared to the period of the signal. Therefore, the total
acoustic pressure can be written as the following sum over harmonics (as
in \cite{du2013fast,soneson2017extending}): \red{
\begin{equation}
    p(\bx, t) = \text{Re} \left\{\sum_{n=1}^{\infty}
        p_n(\bx)e^{-\text{i}n\omega t}\right\}.
    \label{eqn:time_harm_real}
\end{equation}
} Many numerical schemes in the literature consider the time-harmonic
form (\ref{eqn:time_harm_real}), e.g., \cite{campos1999finite,
du2013fast,soneson2017extending, van2015fast}, likely owing to the
distinct advantages of a frequency domain approach:
\begin{itemize}
    \item the challenging task of choosing/developing an efficient time
    stepping scheme can be avoided;
    \item arbitrary frequency-dependent attenuation power laws can be
    easily incorporated (whereas in the time domain one has to contend
    with the fractional Laplacian);
    \item methods such as the boundary element method and volume
    integral equation method, can be applied directly;
    \item the computation regions for higher harmonics can be reduced,
    thereby making simulations more efficient.
\end{itemize}
It is this final point that we study in this article.

The form (\ref{eqn:time_harm_real}) is not well suited for substitution
into (\ref{eqn:west_fractional}) since the right-hand of
(\ref{eqn:west_fractional}) requires the computation of a product. For
such a product, expressions in which real or imaginary parts are
required to be taken lead to cumbersome algebra. Therefore, it is more
straightforward to use the following expression, which is equivalent to
(\ref{eqn:time_harm_real}):
\begin{equation}
    p(\bx,t) = \frac{1}{2}\sum_{n=1}^{\infty}\left(p_n(\bx)e^{-\text{i}n\omega t} +
                                       p_n^*(\bx)e^{\text{i}n\omega t}\right),
    \label{eqn:cc}
\end{equation}
where $^*$ denotes complex conjugation. Substituting (\ref{eqn:cc}) into
(\ref{eqn:west_fractional}) and matching coefficients of $e^{-\text{i}
n\omega t}$ for $n\geq 1$ yields
\begin{equation}
    \nabla^2 p_n
    + k_n^2 p_n = \frac{\beta\omega^2}{2\rho_0 c_0^4}n^2
    \sum_{m=1}^{\infty}p_m(p_{n-m}+2p^*_{m-n}),
    \label{eqn:cascade_1}
\end{equation}
for $n=1,2,\ldots$, where $p_n = 0$ for $n\leq 0$, and the complex
wavenumbers are defined as
\begin{equation}
    k_n = \frac{n\omega}{c_0} + \text{i}\alpha(n\omega).
\end{equation}

The equations (\ref{eqn:cascade_1}) can be further simplified by
neglecting small terms on the right-hand side, which is appropriate in
the weakly nonlinear case, as was considered in, e.g.,
\cite{du2013fast}. Specifically, we neglect all terms $p_ip_j$ and
$p_ip^*_j$ for which $i+j>n$, for $n=1,2,\ldots$, thus giving
\begin{equation}
    \nabla^2 p_n
    + k_n^2 p_n = \frac{\beta\omega^2}{2\rho_0 c_0^4}n^2
    \sum_{m=1}^{n-1}p_m p_{n-m}.
    \label{eqn:cascade_2}
\end{equation}
This is a cascade of inhomogeneous Helmholtz equations in which each
right-hand side is a combination of products of lower harmonics.
Therefore, we can solve the equations sequentially, starting from $n=1$.

We note that the assumption of weak nonlinearity is valid for a range of
\red{FUS} applications for thermal ablation, \red{for example, at a
preclinical stage, when excitation protocols and devices require
characterisation~\cite{kothapalli2018convenient,ries2010real}.} However,
for extremely high focal pressures such as those encountered in
\red{lithotripsy and histotripsy, where pressures much higher than 30MPa
can be seen~\cite{izadifar2017mechanical}}, the terms neglected above
become significant. Therefore, a modified approach would be required,
however we do not consider this case here.

\section{Volume potentials}
\label{sec:volume}

The equations (\ref{eqn:cascade_2}) are each of the general form
\begin{equation}
    \nabla^2 u + k^2 u = f(\bx), \quad \bx\in\mathbb{R}^3.
    \label{eqn:helmholtz}
\end{equation}
Via Green's theorem (see, e.g.,
\cite{colton2013integral,costabel2015spectrum}) it can be seen that the
following volume potential satifies (\ref{eqn:helmholtz}):
\begin{equation}
    u(\bx) = -\int_{\mathbb{R}^3}G_k(\bx,\by)f(\by)\sd \by, \quad \bx\in\mathbb{R}^3,
    \label{eqn:vol_pot}
\end{equation}
where $G_k$ is the fundamental solution to the Helmholtz equation, also
known as Green's function:
\begin{equation}
    G_k(\bx,\by)=\frac{e^{\text{i} k|\bx-\by|}}{4\pi|\bx-\by|},\quad \bx\neq\by.
    \label{eqn:green}
\end{equation}
\red{This function is singular when $\bx=\by$, however integrals of the
function across this singularity may be be evaluated using principal
value techniques or appropriate coordinate transformations, as we shall
see in Section~\ref{ss:compute_potential}.} We note that the integral in
(\ref{eqn:vol_pot}) is over an infinite domain, however in practice we
replace this with a finite domain of integration, $D$. So we have that
\begin{equation}
    u(\bx) = -\int_{D}G_k(\bx,\by)f(\by)\sd \by + \varepsilon(D), \quad \bx\in D\subset\mathbb{R}^3,
    \label{eqn:vol_pot_approx}
\end{equation}
where $\varepsilon(D)$ an error incurred by the introduction of a finite
integration domain. Since the \red{FUS} field is highly focused, a
sensibly chosen finite $D$ will yield a negligibly small error. It is
the purpose of this article to investigate how small $D$ can be made
whilst still yielding accurate approximations to $u$.

For clarity, we write out this integral representation
(\ref{eqn:vol_pot_approx}) for each of the first five harmonics:
\begin{align}
    \label{eqn:harm2} p_2(\bx) &= -\frac{2\beta \omega^2}{\rho_0 c_0^4}
                    \int_{D_2}G_{k_2}(\bx,\by)p_1^2(\by)\sd \by, \\
        \label{eqn:harm3} p_3(\bx) &= -\frac{9\beta \omega^2}{\rho_0 c_0^4}
                    \int_{D_3}G_{k_3}(\bx,\by)p_1(\by) \red{p_2(\by)}\sd \by, \\
        \label{eqn:harm4} p_4(\bx) &= -\frac{8\beta \omega^2}{\rho_0 c_0^4}
        \int_{D_4}G_{k_4}(\bx,\by)(p_2^2(\by)
                                  + 2p_1(\by)p_3(\by))\sd \by, \\
        \label{eqn:harm5} p_5(\bx) &= -\frac{25\beta \omega^2}{\rho_0 c_0^4}
        \int_{D_5}G_{k_5}(\bx,\by)(p_1(\by)p_4(\by)
           + p_2(\by)p_3(\by))\sd \by,
\end{align}
where we have introduced the domains $D_i$, $i=2,3,\ldots$, appropriate
sizes of which are to be determined. In a homogeneous medium, the first
harmonic $p_1$ is merely the incident field generated by the transducer,
which we can compute anywhere in $\mathbb{R}^3$, \red{as is discussed in
Section~\ref{subsec:incident}}. Note finally that the Green's functions
in each of (\ref{eqn:harm2})-(\ref{eqn:harm5}) possess the appropriate
wavenumber for that harmonic, $k_n$, for $n=2, 3, 4, 5$.

\subsection{Efficient computation of the volume potential}
\label{ss:compute_potential}

\red{FUS} problems are reknowned for their challenging high-frequency
nature, with domain sizes up to hundreds of wavelengths in each of the
three dimensions. Therefore, the (singular) integrals in
(\ref{eqn:harm2})--(\ref{eqn:harm5}) are potentially enormously
expensive to approximate.

We choose the integration domains $D_i$ to be cuboidal in shape and to
be discretised into uniform voxel grids $\mathcal{V}(D_i)$ so that the
discrete form of the integral operator becomes a Toeplitz matrix. A
Toeplitz matrix of dimension $N$ has the property that it may be
embedded in a circulant matrix of dimension $2N$, with which a
matrix-vector product can be computed with $\mathcal{O}(2N\log 2N)$
complexity using the fast-Fourier transform (FFT); see, e.g.,
\cite{groth2020accelerating} for more details.

For the voxels in which the Green's function's singularity is located,
we approximate the integral by the integral over a sphere of radius $a$,
chosen such that the sphere's volume is equal to that of the voxel.
\red{This integral can then be transformed into spherical coordinates,
which is convenient because the Jacobian determinant of the
transformation cancels the singularity in the Green's function.} In the
non-singular voxels we approximate the volume potential integral by the
midpoint rule. This gives us the simple quadrature rule
\begin{linenomath}
\begin{align}
    \nonumber\int_{V_j}&G(\bx_i,\by)f(\by)\sd \by \\
    &=\begin{cases}
        \frac{1}{k^2}\{e^{\text{i}ka}(1-\text{i}ka) - 1\}f(x_j), \ &\bx_i\in V_j, \\
        (\delta x)^3 G(\bx_i, \bx_j)f(\bx_j),\ &\bx_i\notin V_j,
    \end{cases}
    \label{eqn:quad}
\end{align}
\end{linenomath}
for $i,j=1,\ldots,N$, where $N$ is the number of voxels in the grid and
$\delta x$ is the side length of each voxel $V_j$. This is reminiscent
of the `discrete dipole approximation' often used in electromagnetic
scattering calculations~\cite{draine1994discrete}. One can choose to opt
for a more sophisticated quadrature rule here, however, this simple
approach suffices for our purposes.

\subsection{\red{FUS} incident field}
\label{subsec:incident}

In a homogeneous medium there is no scattering and therefore the first
harmonic, $p_1$, can be obtained directly from an appropriate model of
the time-harmonic field generated by the \red{FUS} source. In this
paper, we consider bowl-shaped ultrasound transducers, which are
designed to focus acoustic energy at a prescribed location, typically
the centre of curvature of the bowl. More specifically, we consider for
simplicity a single-element bowl-shaped transducer and discretise the
surface using a Rayleigh integral type approach described below. We note
that we are not exploiting symmetries in our approach and therefore more
sophisticated multi-element transducer arrays, such as those in
\cite{gelat2011modelling,gavrilov2000theoretical,kreider2013characterization},
can be incorporated in a straightforward manner.

To discretise the surface of the bowl transducer, we use evenly spaced
points following \cite{deserno2004generate}. At each of the evenly
spaced point we place a monopole source, the expression for which is
given by (\ref{eqn:green}). Summing over the monopole sources gives the
(unnormalised) first harmonic at any $\bx\in\mathbb{R}^3$ as
\begin{equation}
    \tilde{p}(\bx) = \frac{A}{n_{\text{p}}}\sum_{i=1}^{n_{\text{p}}}G_k(\bx,\br_i),
    \label{eqn:incident_field}
\end{equation}
where $n_{\text{p}}$ is the number of points, $\br_i$ are their
locations on the bowl, and $A$ is the total surface area of the bowl.

We note that in (\ref{eqn:incident_field}) no amplitude has been
specified for the monopole sources. We instead choose to normalise the
field to produce a prescribed total radiated power, $\Pi$, which is
obtained by integrating the intensity over a sphere surrounding the
source~\cite{kinsler1999fundamentals}. For a bowl-shaped transducer, the
field is directed, so rather than integrating over a sphere, it suffices
to integrate over a disc covering the open end of the bowl. This can be
written as
\begin{equation}
    \Pi(p) = \frac{1}{2\rho_0 c_0}\int_{0}^{2\pi}\int_{0}^R
    p^2(r,\theta)\sd r \sd \theta,
\end{equation}
where $R$ is the outer radius of the bowl. Thus, the normalised first
harmonic to yield a prescribed radiated power $\Pi_0$ is given as
\begin{equation}
    p_1(\bx) = \sqrt{\frac{\Pi_0}{\Pi(\tilde{p})}}\tilde{p}.
\end{equation}
In our experiments, we take $n_{\text{p}}=4096$, \red{which equates to
approximately ten monopole sources per fundamental wavelength}. Such a
large value was required to avoid undesired interference patterns
between the bowl and focus.

The two bowl transducer geometries considered in this article are taken
from the Sonic Concepts website~\cite{sonic} and are detailed in
Tab.~\ref{tab:transducers}.
\begin{table}[h!]
    \centering
    \begin{tabular}{c  c  c  c}
        \hline\hline
           & $f_0$ (MHz) & $l$ (mm) & $R$ (mm)\\ 
        \hline
        H101 & 1.1 & 63.2 & \red{32} \\ 
        H131 & 1.1 & 35 & \red{16.5} \\ 
        \hline\hline
    \end{tabular}
    \caption{Operating frequencies and geometrical parameters of bowl
    transducers considered, taken from \cite{sonic}. The geometrical
    parameters are the geometric focal length $l$ and the outer radius
    $R$.
    These can also be seen depicted in Fig.~\ref{fig:domain_dimensions}.}
    \label{tab:transducers}
\end{table}
Two propagation media are considered throughout: water and liver. The
acoustic parameters for these are detailed in Tab.~\ref{tab:media}.
\begin{table}[h!]
    \centering
    \begin{tabular}{c  c  c  c  c  c}
        \hline\hline
           & $\rho_0$ (kg/m$^3$) & $c$ (m/s) & $\beta$ & $\alpha_0$ & $\eta$ (dB/m) \\
        \hline
        Water & 1000 & 1480 & 3.5 & 0.2 & 2 \\
        Liver & 1060 & 1590 & 4.4 & 90.0 & 1.1 \\
        Kidney & 1050 & 1570 & 4.7 & 10 & 1 \\
        \hline\hline
    \end{tabular}
    \caption{Relevant medium parameters for water, liver and kidney at
    1MHz~\cite{duck2013physical,azhari2010basics}. The parameters
    $\alpha_0$ and $\eta$ pertain to the absorption power law in
    (\ref{eqn:power_law}).}
    \label{tab:media}
\end{table}

\subsection{How many points per wavelength?}

Before deploying our scheme to compute high-order harmonics, it is
necessary to understand the convergence rate as the mesh is refined,
thereby enabling an adequate resolution to be chosen for later
investigations. In this article, we focus on achieving approximations
with relative $L^2$ errors close to 1\%.

\begin{figure}
    \includegraphics[width=\columnwidth]{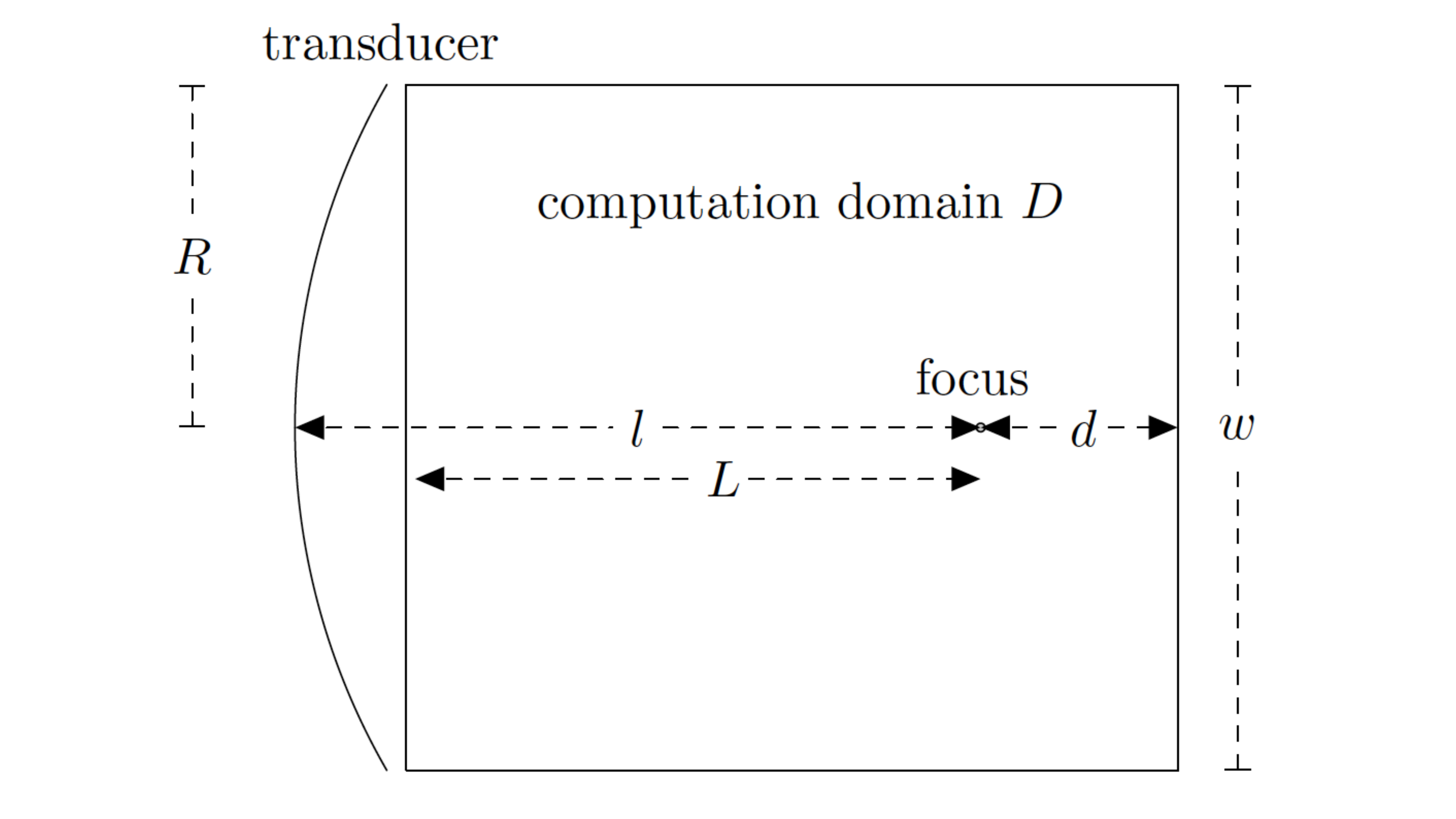}
    \caption{Description of a generic computation domain $D$
    with length $L+d$ and width (and depth) $w$.}
    \label{fig:domain_dimensions}
\end{figure}

\begin{figure}[h!]
    \centering
    \includegraphics[width=\linewidth]{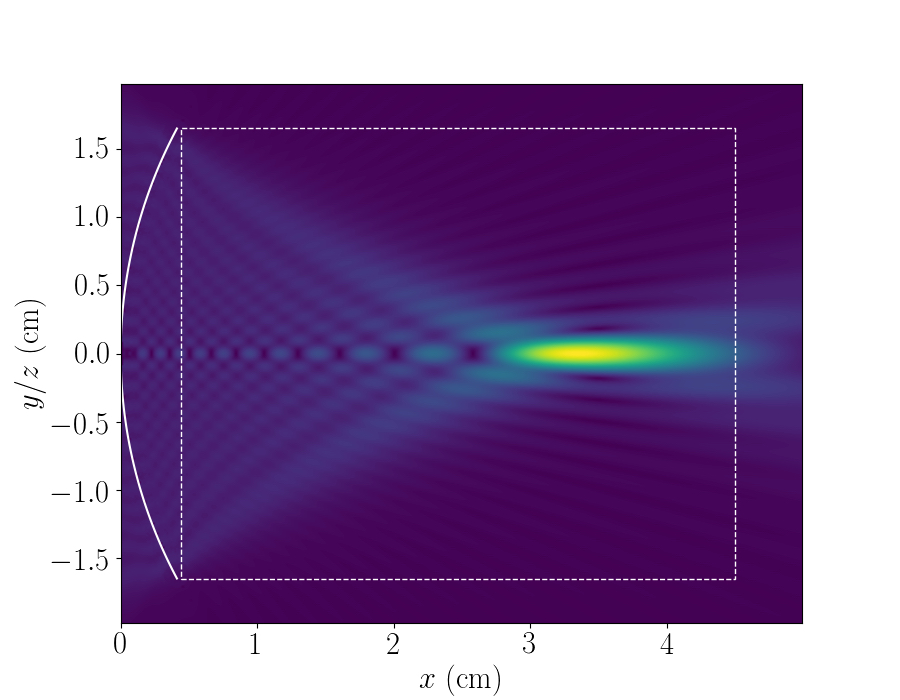}
    \caption{Magnitude of the first harmonic generated by the H131
    transducer in liver. The bowl transducer is represented by the arc
    at the far left. The dashed white line outlines the computation
    domain $D_2$ used to compute the second harmonic shown in
    Fig.~\ref{fig:p2_quad}.}
    \label{fig:H131_p1}
\end{figure}

To determine the convergence rate, we consider computing the second
harmonic generated by the H131 transducer in liver, using equation
(\ref{eqn:harm2}). The configuration of the transducer is shown in
Fig.~\ref{fig:H131_p1}; the bowl transducer is located at the origin and
is directed along the $x$-axis. For the integration domain $D_2$ in
(\ref{eqn:harm2}) we take the white box in Fig.~\ref{fig:H131_p1}. This
box is illustrated with detailed dimension definitions in
Fig.~\ref{fig:domain_dimensions}.
The domain $D_2$ is defined as
\begin{equation}
    D_2 = [l-L, l+d] \times [-R, R] \times [-R, R],
    \label{eqn:domain}
\end{equation}
where $d$ is the distance of interest beyond the focus, $l$ is the
geometric focal length, and $R$ the outer radius (see
Fig.~\ref{fig:domain_dimensions}). The length $L$ is defined as $L =
\sqrt{l^2-R^2} - \epsilon$ with $\epsilon$ chosen as a small
displacement from the bowl \red{in order to avoid the possibility of
$\bx=\br_i$ in (\ref{eqn:incident_field}), for which the monopole
sources are undefined}; we take $\epsilon = 0.1$mm. \red{In this
article, we suppose that the region near the focus is of primary
interest, therefore from a computational point of view, we desire to
shrink the computation domain to be much shorter and narrower than $L$
and $w$, i.e., more localised around the focus, in order to reduce
computational load. The extent to which this shrinking can be done
without losing accuracy is the focus of sections
\ref{sec:results}-\ref{subsec:interpolation}. The values of $L$ and $w$
defined above represent our default computation domain. The value of the
distance $d$ depends on the user's interest in the region beyond the
focus. In the absence of scatterers beyond the focus, the field in this
region after the focus will not affect the intensity at the focus.} In
the present example, the total radiated power of the transducer is set
as $\Pi_0 = 100$W and the operating frequency is $f_0=1.1$MHz.

The volume potential in (\ref{eqn:harm2}) is computed using the method
outlined in Section~\ref{ss:compute_potential} with increasingly refined
voxel meshes. Specifically, we build meshes with voxel dimension $\delta
x = \lambda / (2n_w)$, where $\lambda = c / f_0$ is the fundamental
wavelength and $n_w$ is the `number of voxels per wavelength'. The
factor 2 in the denominator is to account for the fact that we are
computing a field with wavelength $\lambda/2$, i.e., the second
harmonic. The values of $n_w$ considered are from 4 to 20, and a
reference solution, denoted $p_2^r$, is computed using $n_w=35$. For
each value of $n_w$, the relative $L^2$-error of the field $p_2$ along
the $x$-axis is computed; this is defined as
\begin{equation}
    \text{Error} = \frac{||p_2 - p_2^r||}{||p_2^r||}\times 100\%,
    \label{eqn:error_def}
\end{equation}
where $||\cdot||$ denotes the $L^2$-norm along the $x$-axis, i.e.,
\begin{equation}
    ||p_2|| := \left(\int_{l-L}^{l+d}|p_2(x,0,0)|^2\sd x\right)^{1/2}.
\end{equation}
We approximate the integrals in (\ref{eqn:error_def}) using the midpoint
rule with the mesh nodes of the reference solution, $p_2^r$, being used
as the quadrature nodes.

\begin{figure}[h!]
    \centering
    \includegraphics[width=\linewidth]{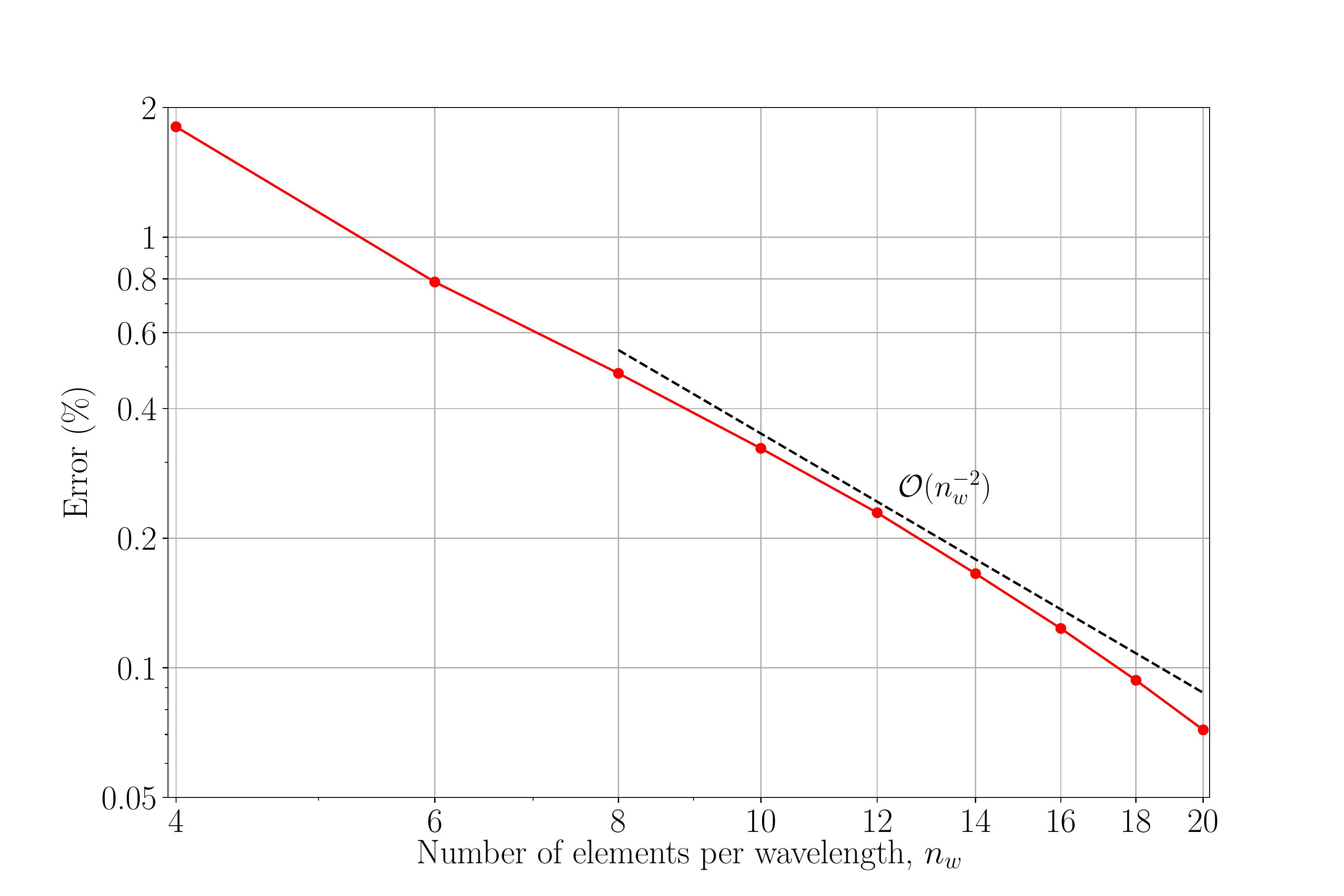}
    \caption{The convergence of the quadrature rule for the computation
    of the second harmonic via (\ref{eqn:harm2}). The convergence rate
    is quadratic in $n_w$ and an error of smaller than 1\% is achieved
    with $n_w>5$.}
    \label{fig:conv_quad}
\end{figure}
\begin{figure}[h!]
    \centering
    \includegraphics[width=\linewidth]{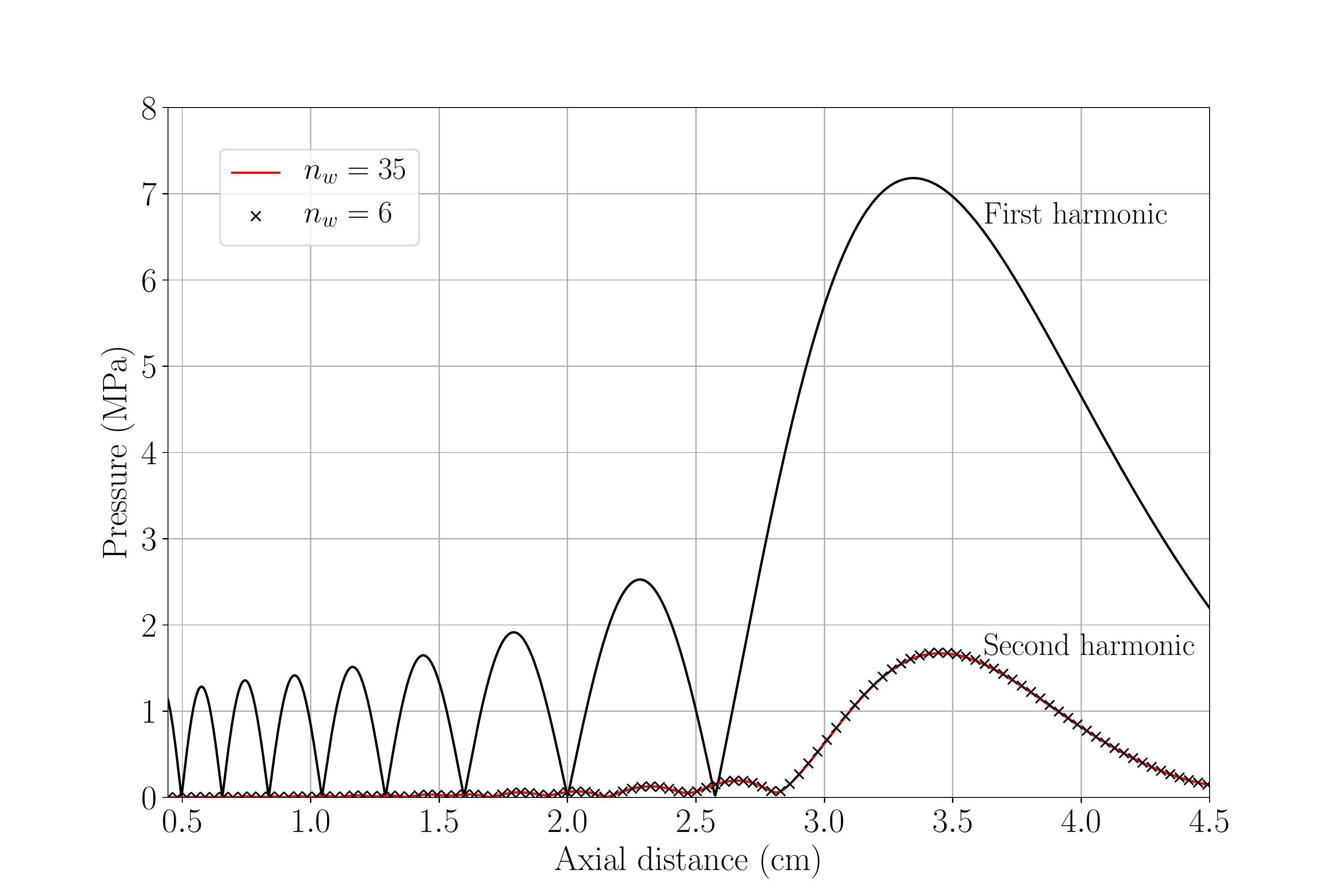}
    \caption{The first and second harmonics along the $x$-axis for the
    H131 transducer in liver. With 6 voxels per wavelength the
    approximation to the second harmonic is indistinguishable from the
    reference solution.}
    \label{fig:p2_quad}
\end{figure}

The convergence results are shown in Figure~\ref{fig:conv_quad}. As is
to be expected from the midpoint rule, quadratic convergence is
obtained. From the graph we can read off that an error smaller than 1\%
is achieved with $n_w>5$. Therefore we choose to take $n_w=6$ for all
harmonics in the experiments in the remainder of the article. The
approximation to the second harmonic with $n_w=6$ is shown in
Fig.~\ref{fig:p2_quad} and can be seen to be indistinguishable from the
reference solution.

\section{Validation of numerical scheme}
\label{sec:validation}

To validate our approach for the computation of higher harmonics via the
evaluation of the volume potentials in
(\ref{eqn:harm2})-(\ref{eqn:harm5}), we present a qualitative
comparison to approximations obtained using \textit{HITU
Simulator}~\cite{HITUwebpage}, \red{which we have chosen because of its
computational efficiency for the axisymmetric problems considered here}.
HITU simulator is an open-source \textsc{Matlab} implementation of the
high-order parabolic approximation to the axisymmetric Westervelt
equation, i.e., the wide-angle Khokhlov-Zabolotkaya-Kuznetsov (WAKZK)
equation. The method is detailed by the author of HITU simulator in
\cite{soneson2017extending}. The assumption of axisymmetry allows the
dimension of the problem to be reduced by one and hence facilitates
rapid simulations.

\begin{figure}[h!]
    \centering
    \includegraphics[width=\linewidth]{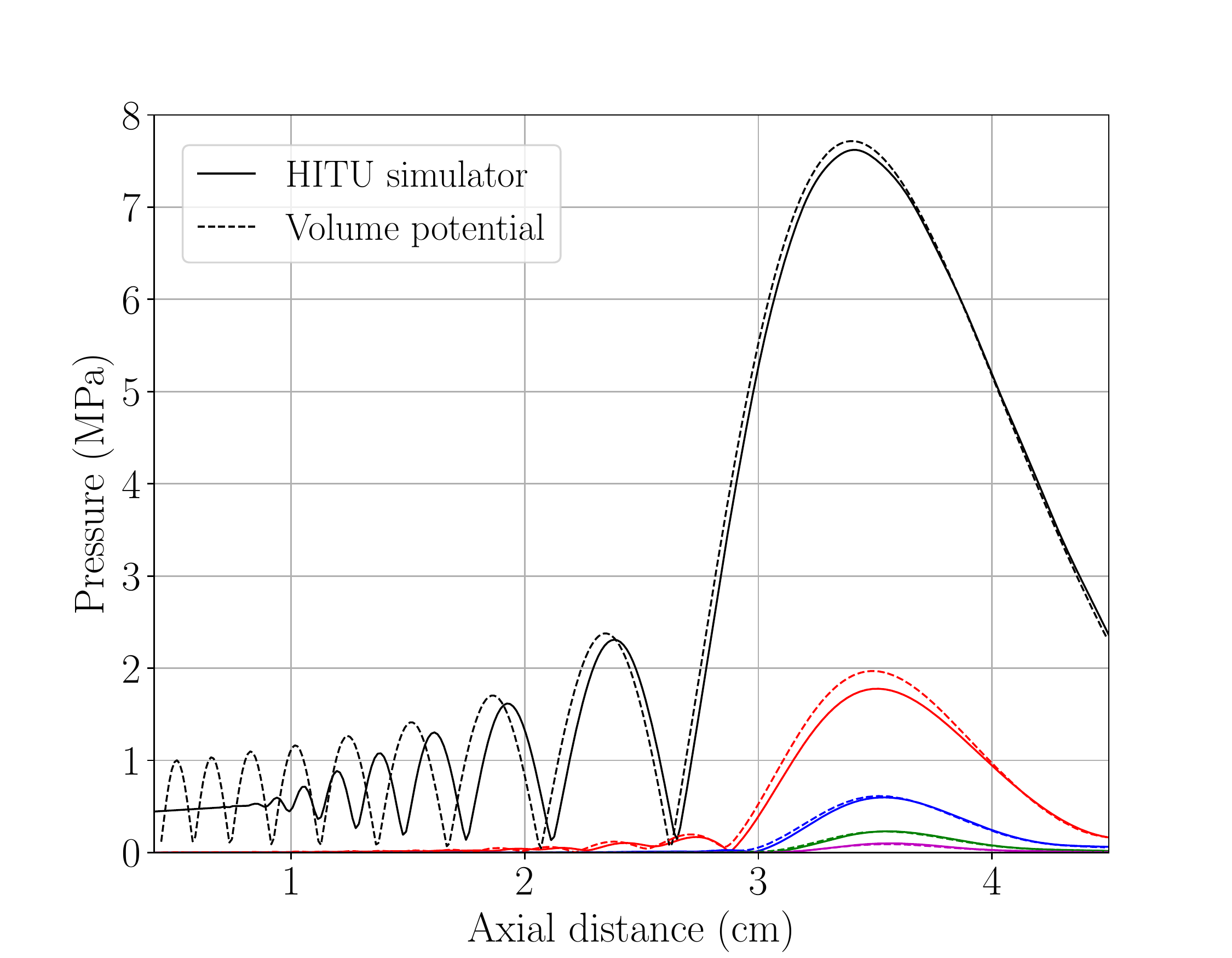}
    \caption{Comparison of the VIE approach with HITU simulator for the
    first five harmonics generated by the H131 transducer operating at a
    power of 50W in water.}
    \label{fig:HITU_comparison_H131}
\end{figure}

\begin{figure}[h!]
    \centering
    \includegraphics[width=\linewidth]{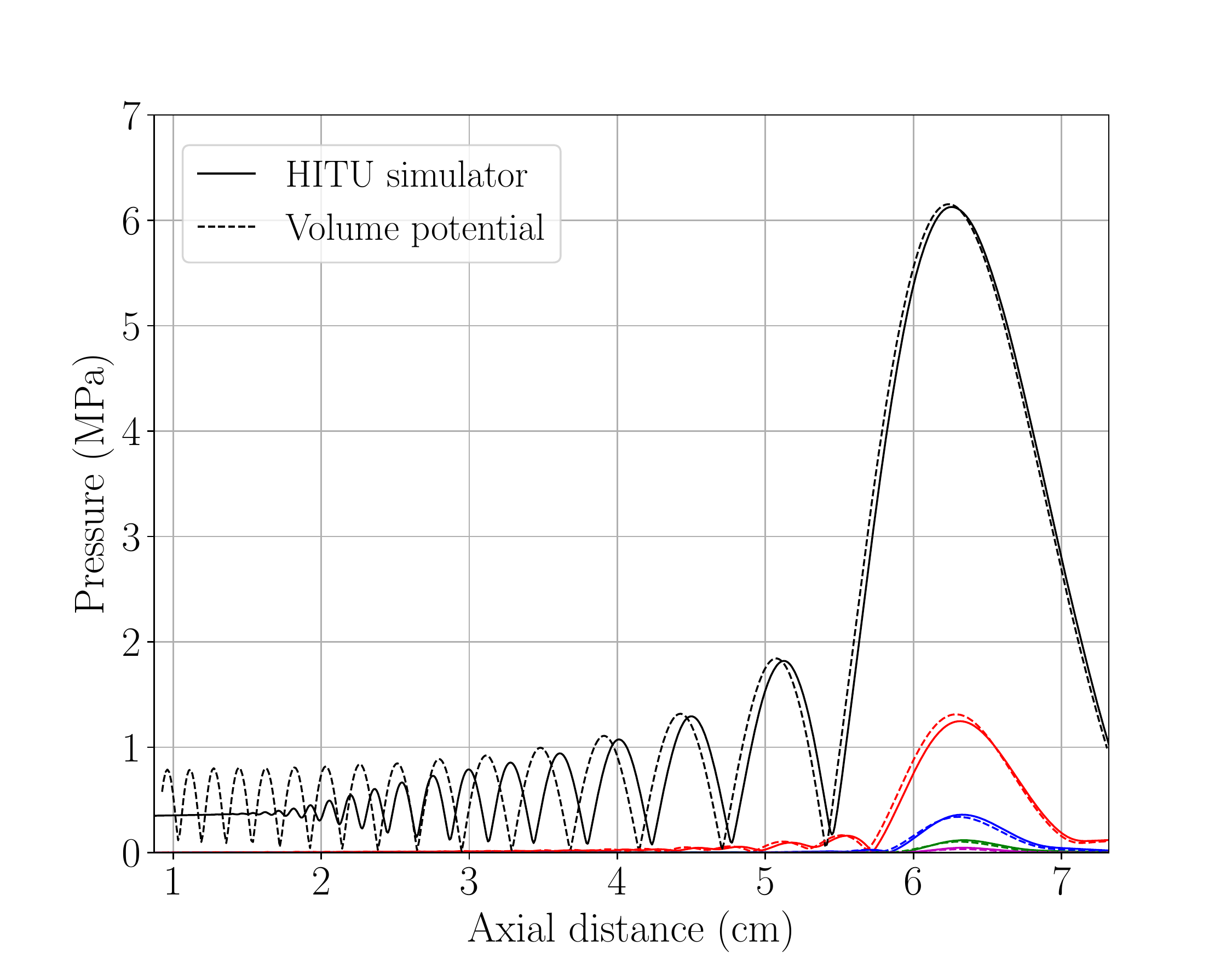}
    \caption{Comparison of the VIE approach with HITU simulator for the
    first five harmonics. H101 transducer at power 100W in liver.}
    \label{fig:HITU_comparison_H101}
\end{figure}

We consider two configurations:
\begin{itemize}
    \item H131 transducer at output power 50W in water;
    \item H101 transducer at output power 100W in liver.
\end{itemize}
\red{These power outputs are chosen since they are relevant to thermal
ablation with FUS (e.g., \cite{kothapalli2018convenient,ries2010real})
and are sufficiently low to ensure that we are in the weakly non-linear
setting.} The first five harmonics along the $x$-axis are shown in
Fig.~\ref{fig:HITU_comparison_H131} and
Fig.~\ref{fig:HITU_comparison_H101}. In both cases we observe good
qualitative agreement with the approximations obtained using HITU
simulator in the region around the focus, whereas towards the transducer
the two methods disagree. This is due to the parabolic assumption made
in the derivation of the WAKZK equation, which is only accurate when
sufficiently far from the transducer. The volume potential method on the
other hand approximates solutions to the Westervelt equation so can be
taken to be accurate in this `near field' region. Indeed, the
computation of the first harmonic in our approach, as outlined in
Section~\ref{subsec:incident}, is equivalent to a Rayleigh integral
method.

\begin{figure*}[ht!]
    \centering
    \includegraphics[width=0.95\linewidth]{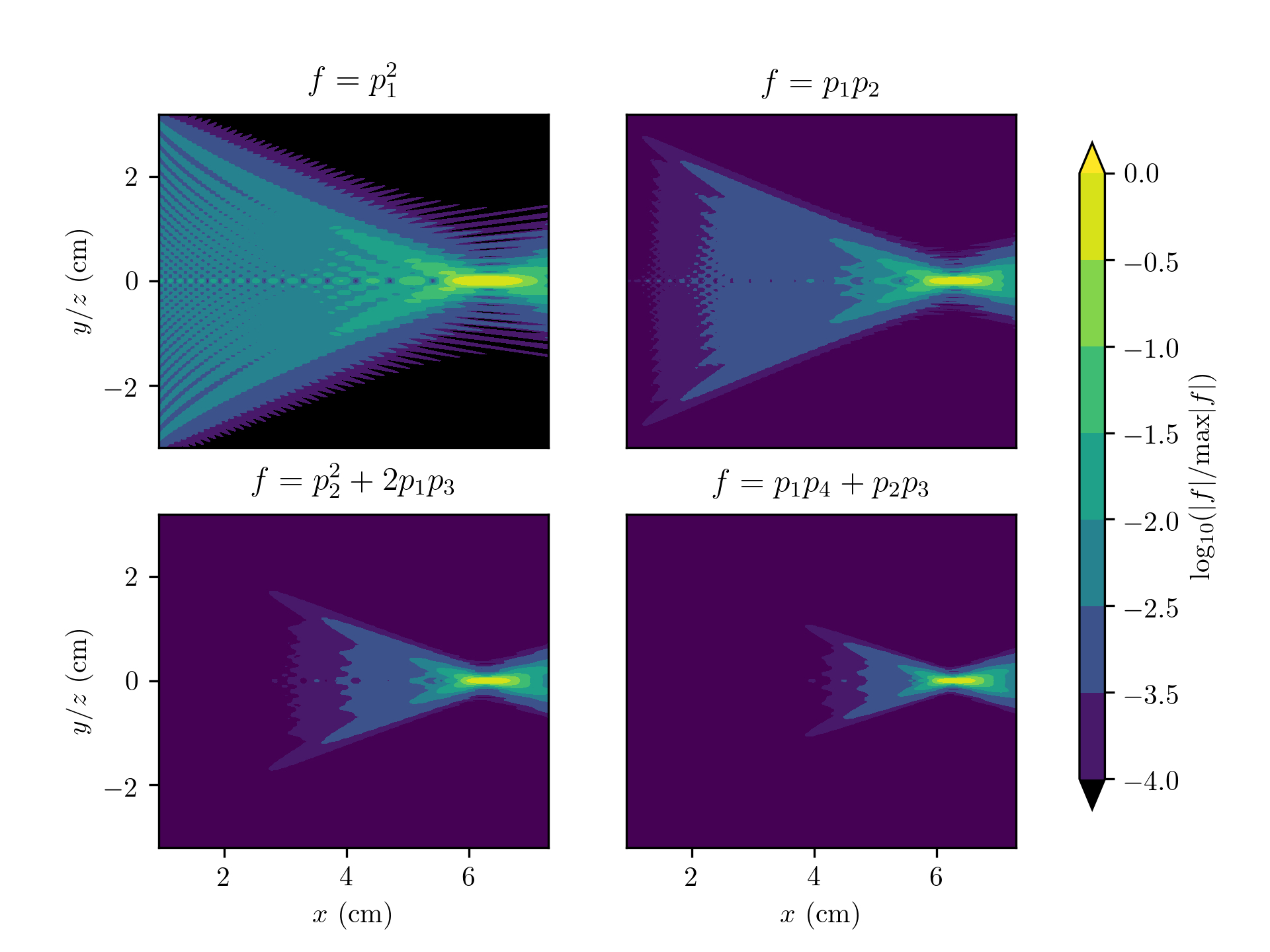}
    \caption{The relative magnitudes of the right-hand side functions
    $f$ as described in equation (\ref{eqn:RHS}) for the H101 transdcuer
    in liver, at 100W. These plots show how the function to be convolved
    with the Green's function becomes more localised as the harmonics
    increase.}
    \label{fig:relative_magnitudes_H101}
\end{figure*}
We observe that the volume potential method predicts a slightly larger
amplitude for the first and second harmonics than HITU simulator,
\red{which may be because the energy flow of the higher harmonics to
lower harmonics is neglected in the derivation of
(\ref{eqn:cascade_2}).} However, the agreement for the third, fourth and
fifth harmonics is almost perfect. The two examples considered have
considerably different attenuation power law parameters, thus the strong
agreement with HITU simulator for both examples demonstrates that the
attenuation is being handled correctly in the volume potential method.

\section{Computation domains for successive harmonics}
\label{sec:results}

In this section we aim to determine the amount by which we can restrict
the computation domain while retaining accurate approximations, and thus
enable acceleration of our simulations.

The main region of interest to practitioners is that around the focus,
since this is where tissue ablation occurs. Ideally one would compute
only on that small region. This, however, leads to the second harmonic
(and thus also the third, fourth, etc.) being poorly approximated since
these harmonics are generated by the accumulation of acoustic energy
over distance. So we seek a balance between accuracy and computational
cost. Here we shall aim to keep the error in each harmonic below 1\%
(relative to the magnitude of the first harmonic) whilst shrinking the
integration regions $D_i$ in (\ref{eqn:harm2})-(\ref{eqn:harm5}) as much
as possible.

It is important to note that (in this homogeneous setting), the field
beyond the focal point has no influence on the field in front of it,
i.e., the waves propagate only in the positive $x$-direction. This means
that our exploration of domain shrinking only applies to the region
between the transducer and the focal point. Beyond the focus we keep the
domain length in the $x$-direction fixed. In most \red{FUS} settings,
the practioners are not interested in the field far beyond the focal
region, therefore the inability to shrink the region beyond the focus
does not greatly affect the gains achieved from shrinking the
computational domain before the focus.

Each of the equations (\ref{eqn:harm2})-(\ref{eqn:harm5}) has the form
\begin{equation}
    p(\bx) = C\int_D G_k(\bx,\by) f(\by) \sd\by,
    \label{eqn:RHS}
\end{equation}
where $f = p_1^2, p_1p_2,\ldots$, $C$ is the appropriate constant, and
$k$ is the appropriate wavenumber. In order to accurately approximate
$p$, the integration domain $D$ must enclose the region where the
integrand is non-negligible. Outside of this region, we can discard the
contributions. The magnitude of the integrand is dictated by the
function $f$, which has the units of intensity. In order to have an idea
of how localised the different functions $f$ are, we plot them for a
particular example in Fig.~\ref{fig:relative_magnitudes_H101}. The setup
considered in the figure is the H101 transducer at 100W in liver. The
figure shows the magnitudes of $f$ scaled by their maximum values (at
the focus) and converted to a log-scale. Consider the top-left image:
this is the $f$ required for the computation of the second harmonic. We
can see that the magnitude is significant all the way back to the
transducer, implying that we must include all this area in the
integration domain $D$. For the remaining images, corresponding to the
third, fourth and fifth harmonics, the functions become increasingly
more localised in both the $x$ and $y/z$ dimensions, suggesting that the
required integration domain can be considerably smaller than that for
the second harmonic.

To investigate this more rigorously, we perform convergence tests for
each harmonic as the relevant integration domain $D$ is restricted. That
is, we take the harmonics generated on the domain in (\ref{eqn:domain})
as the reference solutions $p_l^r$, $l=2,3,4,5$, and then compute the
same harmonics on successively smaller domains and compare the
approximations to the reference. The error of an approximation $p_l$ is
computed along the $x$-axis as
\begin{equation}
    \text{Error} = \frac{||p_l-p_l^r||}
    {||p_1^r||}\times 100\%.
\end{equation}
Note that the harmonic field in the denominator is that of the first
harmonic. This is done so that the error function incorporates the
diminishing size of successive harmonics. For example, if the fifth
harmonic is negligibly small relative to the first harmonic (and so not
worth calculating), the error will reflect this by being very small.

As a measure of the ``localisedness'' of the functions $f$, we use the quantity
plotted in Fig.~\ref{fig:relative_magnitudes_H101}, which we denote as $Q$:
\begin{equation}
    Q(\bx) := \log_{10}\left(\frac{|f(\bx)|}{\text{max}|f(\bx)|}\right).
    \label{eqn:Q}
\end{equation}
In the convergence tests, the integration domain is chosen as the smallest
cuboidal domain $D$ such $Q(\bx)<Q_0$ for $\bx\notin D$, where $Q_0$ is a
given threshold.

\begin{figure}[h!]
    \centering
    \includegraphics[width=\linewidth]{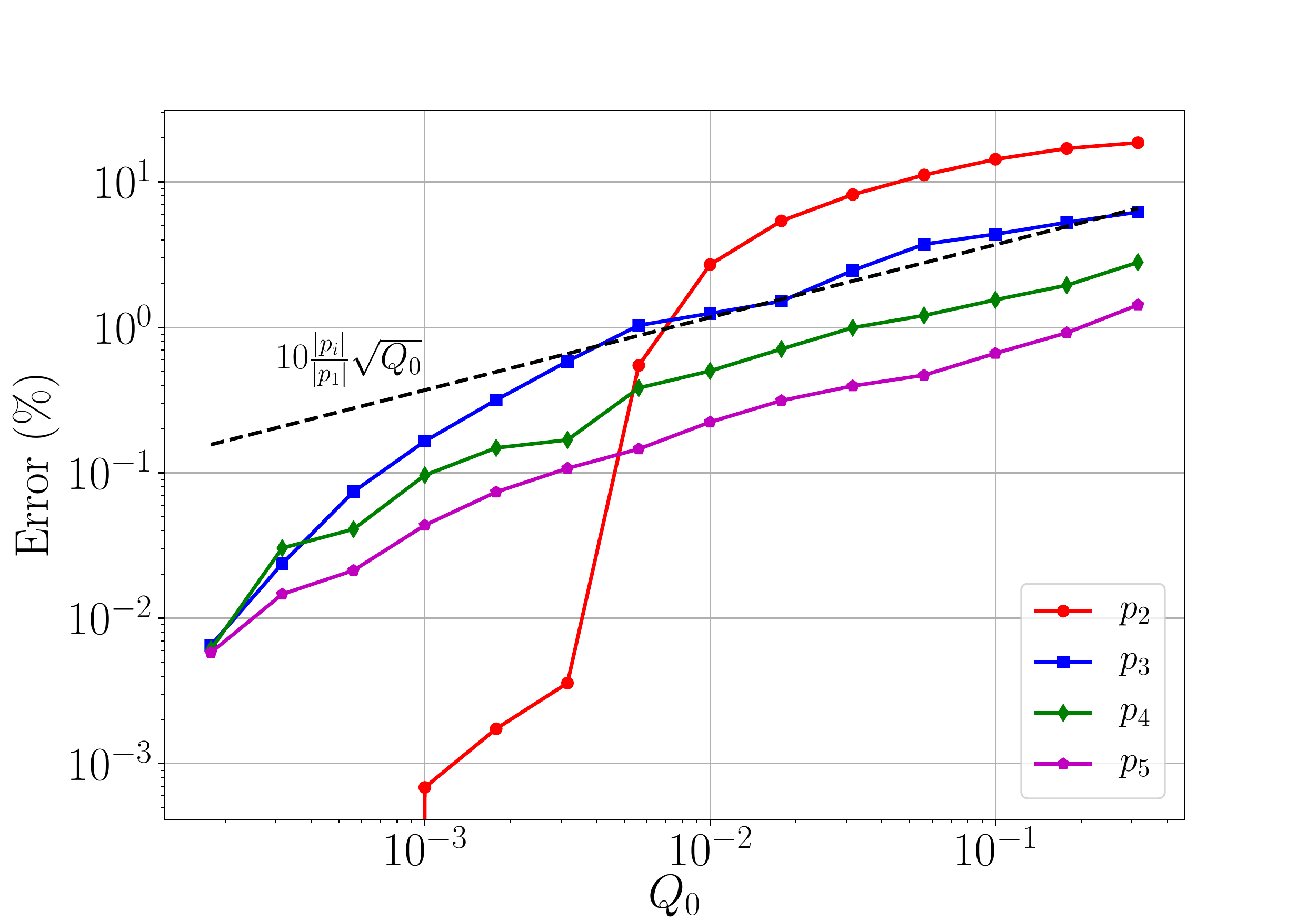}
    \caption{Convergence of the approximations to harmonics $p_i$,
    $i=2,\ldots,5$ as the domain of integration $D_i$ for each is
    adjusted according to the function $Q\geq Q_0$, as defined in
    (\ref{eqn:Q}). The setup considered here is the H131 transducer at a
    power of 100W in water.}
    \label{fig:domain_convergence_H131_water}
\end{figure}

\begin{figure}[h!]
    \centering
    \includegraphics[width=\linewidth]{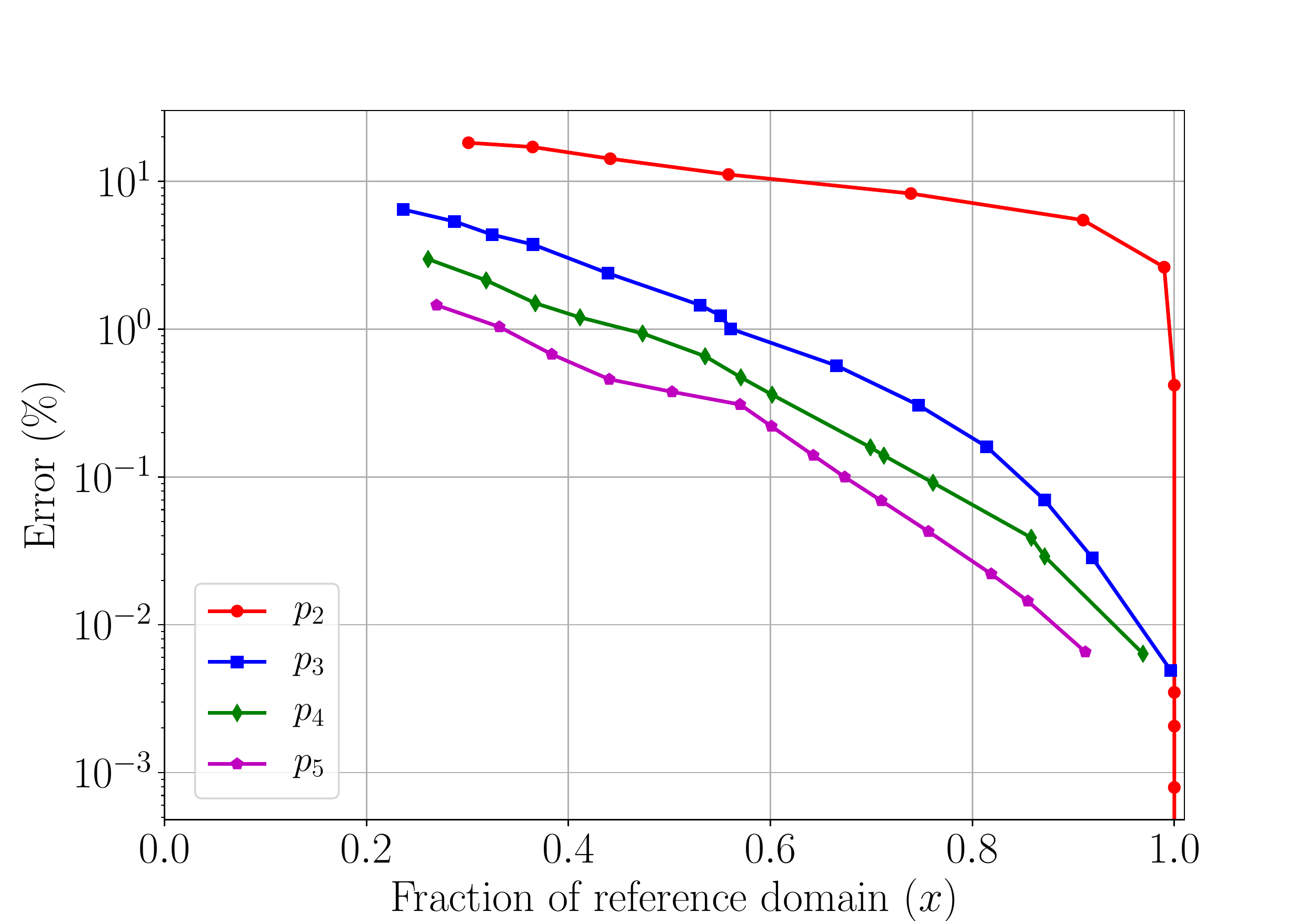}
    \caption{Convergence for H131 transducer at 100W in water, as in
    Fig.~\ref{fig:domain_convergence_H131_water}. The error is plotted
    against the fraction of the total domain in the $x$-direction. This
    demonstrates that, to achieve an error smaller than 1\%, the
    computation domains can be contracted significantly in the
    $x$-direction for harmonics higher than the second.}
    \label{fig:domain_convergence_in_space_H131_water}
\end{figure}

The first configuration we \red{perform the convergence experiment for}
is the H131 transducer in water at an output power of 100W. The
convergence for each harmonic is shown in
Fig.~\ref{fig:domain_convergence_H131_water}. We notice that the
convergence of the second harmonic drops suddenly once the error dips
below 1\% -- this is because the computation domain is close to the size
of the reference domain by this point. This is illustrated more clearly
in Fig.~\ref{fig:domain_convergence_in_space_H131_water} where the same
data as in Fig.~\ref{fig:domain_convergence_H131_water} is shown but now
plotted against the size of the computation domain as a fraction of the
reference domain, rather than against $Q_0$. \red{By `fraction' of the
domain, we mean the scaling factor such that the length $L'$ and width
$w'$ of the shrunken domain are given by
\[
    L' = \text{fraction}_x\times L,\quad w' = \text{fraction}_{y,z}\times w.
\]
}

Fig.~\ref{fig:domain_convergence_in_space_H131_water} shows that, to
achieve less than 1\% error in $p_2$, the computation domain must extend
all the way to the transducer.
For the higher harmonics, however, a different trend is evident. In
Fig.~\ref{fig:domain_convergence_H131_water} we observe that the error
curves for $p_3,p_4,p_5$ each have the approximate behaviour
\begin{equation}
    \text{Error}(p_i) \approx 10\frac{||p_i||}{||p_1||}\sqrt{Q_0},\quad i=3,4,5,
    \label{eqn:trend}
\end{equation}
where $||\cdot||$ represents the $L^2$-norm. Although an interesting
observation, the utility of the relationship (\ref{eqn:trend}) for
dictating an appropriate computation domain for $p_i$ is not immediately
apparent, since it requires the computation of $||p_i||$ before $p_i$
has been computed. Determining an \emph{a priori} approximation for
$||p_i||$ to make use of (\ref{eqn:trend}) would be a useful endeavour,
however we do not undertake such a task in this article. Rather we seek
to develop an approximate rule of thumb for choosing sensibly-sized
computation domains for the harmonics. To this end, it is more
straightforward to consider the convergence of the approximations in
terms of physical distance, as is done in
Fig.~\ref{fig:domain_convergence_in_space_H131_water} and
Tab.~\ref{tab:convergence_domain}.

\begin{figure}[h!]
    \centering
    \includegraphics[width=\linewidth]{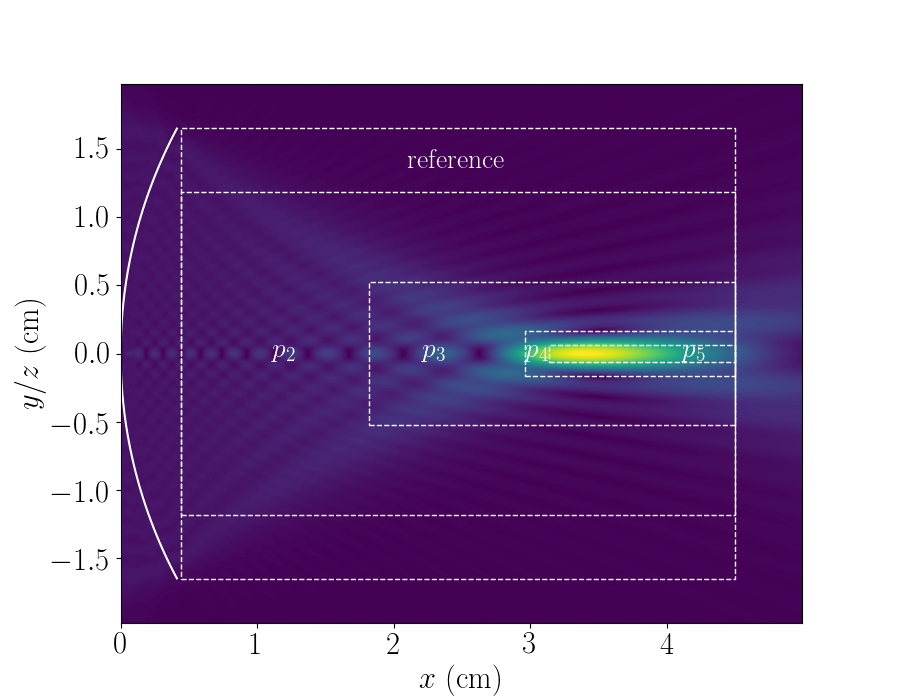}
    \caption{Nested domains for the computation of successive harmonics
    while keeping relative error below 1\%, for the H131 transducer
    operating at 100W in water. The domain used to compute the reference
    solution is in (\ref{eqn:domain}).}
    \label{fig:H131_water_subdomains}
\end{figure}

By looking at the 1\% error line in
Fig.~\ref{fig:domain_convergence_in_space_H131_water}, it is possible to
read off the size of the domains as fractions of the reference domain.
The precise values are reported in Tab.~\ref{tab:convergence_domain}. A
nested series of domains constructed according to this specification is
shown in Fig.~\ref{fig:H131_water_subdomains}. Let us elaborate on the
potential computational gain achieved using these nested domains
compared to computing all harmonics on the reference domain.

The reference domain for the H131 transducer in water has dimensions
$[4.1\text{cm},3.3\text{cm},3.3\text{cm}]$ and is discretised into
voxels of dimension $\delta x = \lambda_5/6\approx \red{45.1\mu}$m,
i.e., fine enough to resolve the highest harmonic of interest. This mesh
has $901\cdot 732\cdot 732\approx 4.8\times 10^8$ voxels for the
computation of all the harmonics. Using the nested domains built
according to the specifications in Tab.~\ref{tab:convergence_domain}, we
obtain a mesh for $p_2$ with dimensions $[4.1, 2.4, 2.4]$cm, which is
discretised into voxels of dimension $\delta x = \lambda_2/6\approx
\red{113\mu}$m. This mesh has $360\cdot 211\cdot 211\approx 1.6\times
10^7$ voxels, which is a factor of 30 smaller than the reference mesh.
\red{The meshes for the computation of $p_3, p_4$ and $p_5$ have
$1.22\times10^7, 6.1\times10^6$ and $4.5\times10^6$ voxels,
respectively. Thus the total number of voxels required for all harmonics
is $4\times4.8\times 10^8=19\times 10^8$ for the single mesh approach
and $3.9\times 10^7$ for the nested meshing (summing over voxels in each
of the four meshes). Hence we achieve close to a factor 50 reduction in
computational load.}

\begin{table}[h!]
    \centering
    \begin{tabular}{c | c  c  c  c  c}
        \hline\hline
             &     & $p_2$ & $p_3$ & $p_4$ & $p_5$ \\
        \hline
        H131 water & $x$   & 1 & 0.67 & 0.47 & 0.38 \\
        100W & $y/z$ & 0.74 & 0.39 & 0.20 & 0.13 \\
        \hline
        H131 water & $x$   & 1 & 0.67 & 0.58 & 0.52 \\
        150W & $y/z$ & 0.71 & 0.40 & 0.39 & 0.21 \\
        \hline
        H131 liver & $x$   & 1 & 0.56 & 0.25 & 0.26 \\
        100W & $y/z$ & 0.90 & 0.20 & 0.05 & 0.06 \\
        \hline
        Rule of & $x$  &  1 & 0.75 & 0.65 & 0.61 \\
        thumb   & $y/z$ & 1 & 0.67 & 0.5 & 0.4
    \end{tabular}
    \caption{Sizes of domains required to achieve less than 1\% error
    for each harmonic as fractions of the reference domain
    (\ref{eqn:domain}). Since all domains are boxes, the fractions of
    the distances in the $x$ and $y/z$ directions are provided.}
    \label{tab:convergence_domain}
\end{table}

This improvement is impressive, however it is easily seen in
Tab.~\ref{tab:convergence_domain} that this particular domain scaling
does not apply to all transducer and material configurations. Therefore,
we decide upon a simple rule of thumb for designing the separate
computational domains that leads to domains greater than or equal to
those given in Tab.~\ref{tab:convergence_domain}. Therefore, the errors
incurred in using these domains will be even smaller than those obtained
when using the ideal domains specified in the table.

A desirable rule for domain sizes would be the following:
\begin{quoting}
    Let $[L+d, w, w]$ be the dimensions of domain $D_2$ (see
    Fig.~\ref{fig:domain_dimensions}) for harmonic $p_2$, with
    corresponding wavelength $\lambda_2$. Then choose the dimensions of
    domain $D_{i}$, $i>2$, as $[\frac{\lambda_i}{\lambda_2}L+d,
    \frac{\lambda_i}{\lambda_2}w, \frac{\lambda_i}{\lambda_2}w]$. That
    is, the domain is scaled according to the wavelength of the harmonic
    being considered.
\end{quoting}
With this rule for creating the domains, we observe that the number of
voxels in each mesh will almost be the same, save for a slight increase
due to the fact that the distance $d$ is not being scaled. This amounts
to a reduction in the overall number of DOF by a factor of approximately
$(n/2)^3$ (in fact, slightly lower than this due to the unscaled portion
of length $d$), where $n$ is the number of harmonics being computed.

In the next section, we outline an algorithm for evaluating the volume
potentials over a set of nested domains constructed according to the
rule of thumb proposed above.

\begin{table*}[ht!]
    \centering
    \begin{tabular}{c  c  c  c  c c}
        \hline\hline
           &  N$^\circ$ voxels & Meshing & Interpolation & Evaluate $G_{k_i}$ & Compute $p_i$\\
        \hline
        $p_2$ & $1.86\times10^8$ & 26.1s & 24m30s & 4m54s & 3m48s \\
        $p_3$ & $2.01\times10^8$ & 22.2s & 3m9s & 5m14s & 4m9s\\
        $p_4$ & $2.15\times10^8$ & 23.9s & 5m50s & 5m27s   & 4m8s   \\
        $p_5$ & $2.30\times10^8$ & 25.5s & 7m12s & 6m7s   & 7m5s \\
        \hline\hline
    \end{tabular}
    \caption{Performance details for the volume potential approach on
    nested meshes with a resolution of six voxels per wavelength (for
    each harmonic). The configuration considered is the H101 transducer
    operating at 100W in water. Total time taken 1h23m11s.}
    \label{tab:performance}
\end{table*}

\section{Volume potentials on nested meshes}
\label{subsec:interpolation}
To demonstrate the effectiveness of the volume potential evaluation on nested meshes,
we present some final results detailing the computational performance of this
approach. First we outline the algorithm for computing the first $n$ harmonics:
\begin{algorithmic}
    \For {$i=2\rightarrow n$}
    \begin{enumerate}
        \item Create domain $D_i$ for $p_i$ and voxel mesh $\mathcal{V}(D_i)$
        \item Assemble components for integration (\ref{eqn:quad}):
            \begin{enumerate}
                \item Evaluate/interpolate $p_{i-1},p_{i-2},\ldots,p_1$
                at voxel centres in $\mathcal{V}(D_i)$
                \item Evaluate integral of Green's function $G_{k_i}$
                over each voxel
            \end{enumerate}
        \item Compute $p_i$ via appropriate equation in
        (\ref{eqn:harm2})--(\ref{eqn:harm5})
    \end{enumerate}
    \EndFor
\end{algorithmic}
Note that in the above algorithm, for $p_2$, the pressure field $p_1$ is
evaluated over the voxel mesh $\mathcal{V}(D_2)$ as described in
Section~\ref{subsec:incident}. Whereas for later harmonics we perform
interpolation of the earlier harmonics down onto the new mesh. In this
work we have used linear interpolation, however if higher accuracy is
required, we recommend quadratic interpolation (albeit at a higher
computational cost). The second and third step each contain applications
of the FFT: for the circulant embedding of the Green's function in step
2 (see, e.g., \cite{groth2020accelerating}) and to perform the
convolution required in the quadrature rule (\ref{eqn:quad}) in step 3.
These FFTs are performed using the Python wrapper `pyfftw' to the FFTW
library~\cite{FFTW05}.

\begin{figure}[h!]
    \centering
    \includegraphics[width=\linewidth]{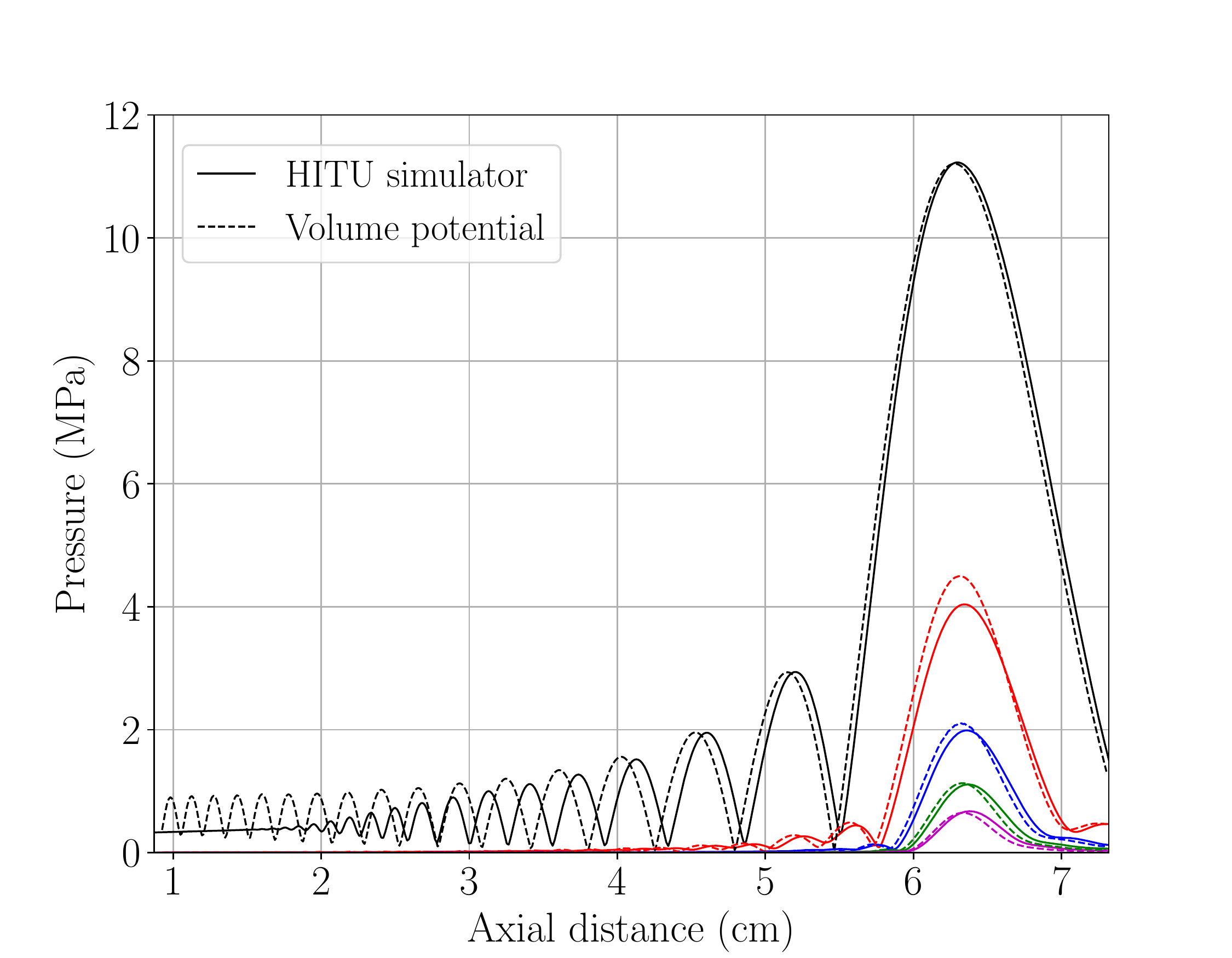}
    \caption{The on-axis absolute pressure field for the first five
    harmonics generated by the H101 transducer at 100W in water, as
    computed using the volume potential approach on nested domains,
    designed to keep the relative error below 1\%. The approximation
    obtained using the HITU simulator is provided for comparison.}
    \label{fig:HITU_comparison_H101_water}
\end{figure}

As an example, we consider the H101 transducer operating at 100W in
water run on a workstation with two sockets, each containing a 14-core
Xeon E5-2690 v4 CPU, each supporting hyper-threading with two threads,
and hence a total number of 56 threads. The total amount of RAM
available on this machine is approximately 270GB, which is ample for the
problems considered here. The first five harmonics along the $x$-axis
are shown in Fig.~\ref{fig:HITU_comparison_H101_water}, along with the
approximation obtained with HITU simulator; again we see that the volume
potential approach predicts a larger peak in the second harmonic but in
general a good qualitative agreement is observed.

To get a feel for the performance of our approach (of evaluating the
volume potentials on nested meshes designed according to our proposed
rule of thumb), the cost of each step in the above algorithm is detailed
in Tab.~\ref{tab:performance}. The largest mesh has 230 million voxels,
as compared to 2.9 billion voxels without nested meshing. This
represents a large saving, in both time and memory. The most expensive
step reported in Tab.~\ref{tab:performance} is the evaluation of the
first harmonic $p_1$ on $\mathcal{V}(D_2)$. This process was described
in Section~\ref{subsec:incident} and has complexity $\mathcal{O}(n_e
N)$, where $n_e$ is number of points used to discretise the surface of
the transducer and $N$ is the number of voxels. Since we take
$n_e=4096$, this is a rather expensive procedure, even when parallelised
over 56 threads. Therefore a transducer model using far fewer elements
and/or a machine with more threads will lead to a large reduction in
computation time for this step. Linear interpolation is computed
efficiently using the scipy~\cite{jones2001scipy}
\verb!RegularGridInterpolator! command. The evaluation of $G_{k_i}$
consists of a parallelised loop over the $N$ voxels and then an FFT of a
three-dimensional complex-valued array of size $(2N_x,2N_y,2N_z)$, where
$N_x,N_y,N_z$ are the number of voxels in each dimension. Then the
computation of $p_i$ consists of one forward and one inverse FFT of an
array of size $(2N_x,2N_y,2N_z)$. The FFTs take advantage of
multithreading, therefore can be easily accelerated through the use of
more cores. Furthermore, it seems likely that optimising the FFT routine
for the particular setting and using the FFTW C++ library directly will
lead to further acceleration. Nevertheless, the current implementation
yields fast and accurate predictions of the first five harmonics (in the
cases considered).

\red{
\section{Extension to inhomogeneous domains}
We briefly consider the extension of this approach to the case of an
inhomogeneous domain. That is, we suppose that the wavespeed, $c(\bx)$,
and non-linearity parameter, $\beta(\bx)$, are now spatially varying.
Further, we assume that the density is close to constant, i.e.,
$\rho(\bx)\approx\rho_0$ (we comment on large density contrasts at the
end of this section). The spatial variation of $c$ and $\beta$ lead to
backscattering of the field generated by the transducer, and thus,
rather than computing the harmonics via direct evaluation of volume
potentials as before, we must now in addition solve \textit{volume
integral equations} (VIEs) account for the scattering effects.

Let the first harmonic of the \textit{incident field} generated by the
transducer be denoted as $p^{\text{inc}}$. Then VIEs for the first five
harmonics are (see \cite{costabel2015spectrum}):
\begin{align}
    \label{eqn:harm1_inhomo} p_1(\bx) -\int_{D_1}G_{k_1}(\bx, \by)(k_1^2(\by)-\overline{k}_1^2)p_1(\by)\sd\by
    &= p^{\text{inc}}(\bx), \\
    \label{eqn:harm2_inhomo} p_2(\bx) -\int_{D_2}G_{k_2}(\bx, \by)(k_2^2(\by)-\overline{k}_2^2)p_2(\by)\sd\by
    &= -\frac{2\beta(\bx) \omega^2}{\rho_0 c(\bx)^4}
                    \int_{D_2}G_{k_2}(\bx,\by)p_1^2(\by)\sd \by, \\
        \label{eqn:harm3_inhomo} p_2(\bx) -\int_{D_3}G_{k_3}(\bx, \by)(k_3^2(\by)-\overline{k}_3^2)p_3(\by)\sd\by &=
        -\frac{9\beta(\bx) \omega^2}{\rho_0 c(\bx)^4}
                    \int_{D_3}G_{k_3}(\bx,\by)p_1(\by) p_2(\by)\sd \by, \\
        \label{eqn:harm4_inhom} p_4(\bx) -\int_{D_4}G_{k_4}(\bx, \by)(k_4^2(\by)-\overline{k}_4^2)p_4(\by)\sd\by &
        = -\frac{8\beta(\bx) \omega^2}{\rho_0 c(\bx)^4}
        \int_{D_4}G_{k_4}(\bx,\by)(p_2^2(\by)
                                  + 2p_1(\by)p_3(\by))\sd \by, \\
        \label{eqn:harm5_inhomo} p_5(\bx) -\int_{D_5}G_{k_5}(\bx, \by)(k_5^2(\by)-\overline{k}_5^2)p_5(\by)\sd\by &=
        -\frac{25\beta(\bx) \omega^2}{\rho_0 c(\bx)^4}
        \int_{D_5}G_{k_5}(\bx,\by)(p_1(\by)p_4(\by)
           + p_2(\by)p_3(\by))\sd \by,
\end{align}
where $\overline{k}_i$ is the wavenumber of the background medium for
harmonic $i$, and $k_i(\bx)$ is the variable wavenumber. Note that the
integrals on the left-hand sides have non-zero contributions only where
$k_i(\bx)\neq \overline{k}_i$, i.e., where the wavenumber differs from
that of the background medium. If $k_i(\bx)\equiv \overline{k}_i$, then
we are in the homogeneous case considered before.

To validate the rule of thumb for an inhomogeneous medium, we consider a
2cm layer of kidney tissue surrounded by water. The layer is centred at
the focus of the transducer. As the tissue properties for water and
kidney, we use those given in Table~\ref{tab:media}, except that we
assume the density of kidney to be equal to that of water, to coincide
with our constant density assumption. We consider the H131 transducer
operating at 50W. A comparison with the HITU simulator is presented in
Figure~\ref{fig:HITU_comparison_slab} where we observe good agreement in
terms of focus location as well as magnitudes of the separate harmonics.
Since HITU is a one-way solver, it does not approximate the
backscattering, whereas our full-wave solver does. The backscattering
can be observed as the ripples in the VIE curves.

\begin{figure}[h!]
    \centering
    \includegraphics[width=\linewidth]{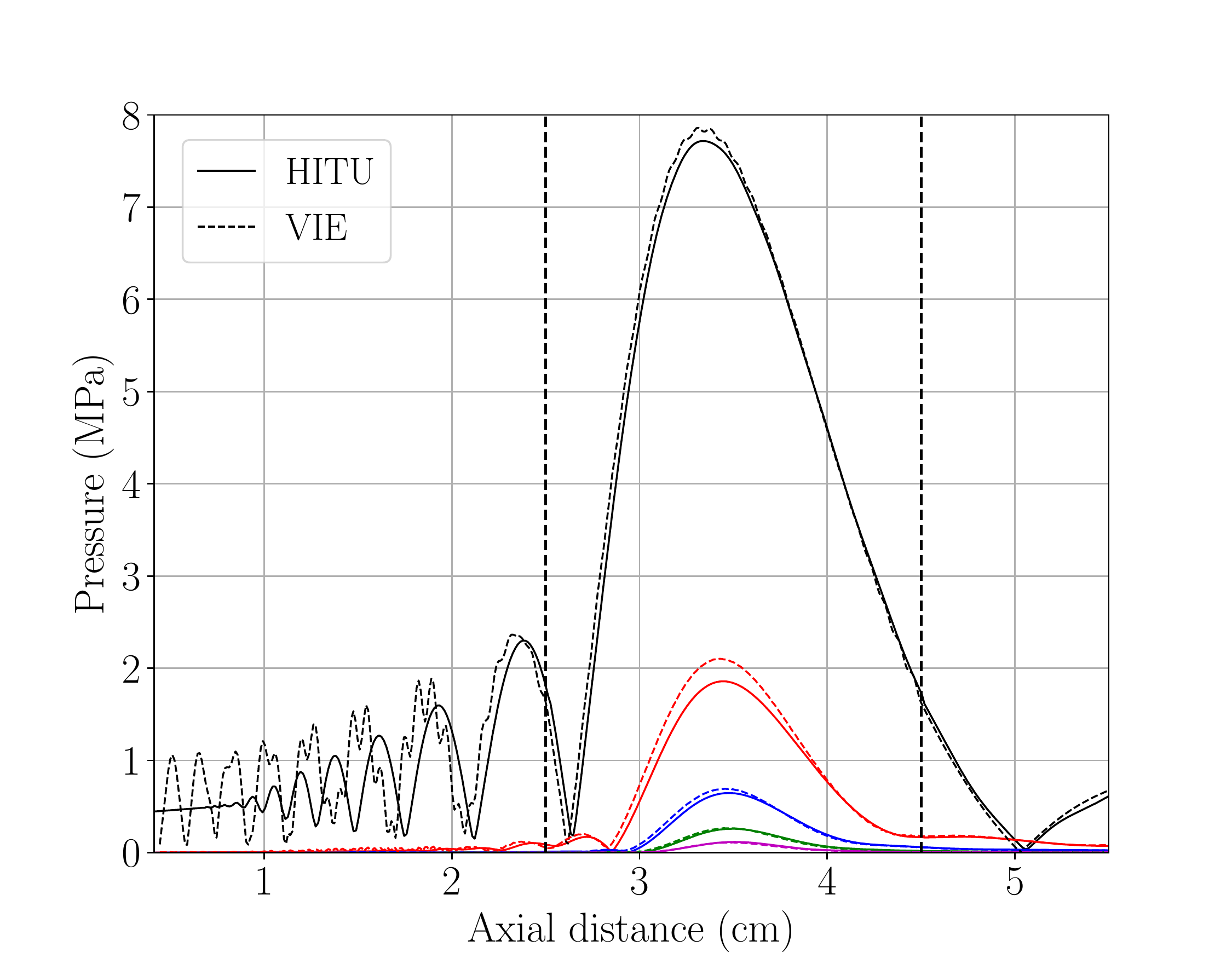}
    \caption{The on-axis absolute pressure field for the first five
    harmonics generated by the H131 transducer at 50W in water with a
    2cm layer of kidney material centred at 3.5cm. The vertical dashed
    lines demarcate the kidney layer.}
    \label{fig:HITU_comparison_slab}
\end{figure}

We note that in the above we assumed that $\rho(\bx)\approx\rho_0$
throughout the inhomogeneous domain. This was done in order to derive
convenient volume integral equations, which can be solved in an
efficient manner. For strong density contrasts, the VIEs
(\ref{eqn:harm1_inhomo})-(\ref{eqn:harm5_inhomo}) must be augmented with
\textit{boundary integrals}, as discussed in
\cite{costabel2015spectrum}. This complicated could be resolved by a
coupling to an established boundary element code, such as
\cite{van2015fast}, but this is left to future work.}

\section{\label{sec:conclusion}Conclusion}

In this paper we have set out
to reduce the computational burden of numerical schemes for \red{FUS}
simulations through the construction of an efficient and simple meshing
strategy. This strategy can be employed with those numerical schemes
that \red{seek to approximate the Westervelt equation on a single
non-uniform mesh}, and those that solve for each harmonic on separate
meshes, such as the frequency-domain volume potential approach proposed
here.

The strategy exploits the increasingly localised nature of the higher
harmonics around the transducer's focal region so that the degrees of
freedom in the mesh can be more efficiently distributed. If we were
considering a single non-uniform mesh approach \red{in which we
approximate the full Westervelt equation}, this mesh would become
increasingly more refined toward the focus, since this is where the
higher harmonics are present. In the frequency-domain setting we
considered, this leads to a nested series of meshes, as was discussed in
detail in this article.

In the frequency domain, the Westervelt equation can be rewritten as a
series of inhomogeneous Helmholtz equations. When the propagation medium
is taken to be homogeneous, these Helmholtz equations can be solved
exactly by volume potentials, which may be efficiently evaluated using
the quadrature method proposed in Section~\ref{sec:volume}. This novel
application of this approach allows us to explore efficiently the
convergence of each harmonic as the respective computation domain was
changed in size. Thus enabling us to determine the smallest domains we
could use in order to achieve an error of less than 1\%.

We showed that the accurate approximation of the second harmonic
requires a computation domain that extends from the focus all the way to
the transducer, since the first harmonic is not sufficiently localised
near the focus to allow a smaller domain to be employed. The third
harmonic and above, however, can be approximated accurately on
considerably reduced domains. We found that scaling the computation
domain's width and height relative to the wavelength under consideration
allowed for accurate approximations for the first five harmonics for the
\red{FUS} configurations considered here. This leads to a reduction in
the number of degrees of freedom of approximately $(n/2)^3$, where $n$
is the number of harmonics being computed.

\red{Finally, we demonstrated how this approach generalises, via the
 introduction of volume integral operators, to inhomogeneous media with
 low density variation. The application to inhomogeneous media with
 large density contrast, such as between water and bone, requires the
 introduction of further boundary integral operators, with is left to
 future work.}

\red{To conclude, we briefly comment on the generalisation of the `rule
of thumb' to other transducer configurations and frequencies. In the
present article, two different transducers were considered, both at
1.1MHz, and propagating within three different media: water, liver,
kidney. All the examples considered produced peak amplitudes of lower
than 15MPa, at which the weakly nonlinear assumption used in the
derivation of the cascade of Helmholtz equations is accurate. We believe
that for different focused transducer configurations and frequencies,
our proposed rule of thumb is accurate, provided the field can still be
categorised as weakly nonlinear. For highly nonlinear fields, a further
study would be required to test the `rule of thumb', and also a
modification of our volume potential approach required to allow for the
transfer of energy from higher to lower harmonics.}

A Python implementation of this work is freely available at
\begin{center}
    \verb!github.com/samuelpgroth/vines!.
\end{center}

\section*{Acknowledgments}
This work was supported by a grant entitled ``Optimising patient
specific treatment plans for ultrasound ablative therapies in the
abdomen (OptimUS)'' from the EPSRC (EP/P013309/1 to Cambridge,
EP/P012434/1 to UCL).

\bibliographystyle{acm}
\bibliography{Cambridge_bib}

\end{document}